\documentclass{amsart}   %{tran-l}
\usepackage{geometry}
\usepackage{enumerate}
\usepackage{verbatim}

\usepackage{amscd}     % commutative diagrams -- see page 293
\usepackage{amsthm}   % theorem-like definitions
\usepackage{amssymb}
\usepackage{amsmath}
\usepackage{color}
\usepackage[all]{xy}
\copyrightinfo{2009}{American Mathematical Society}

\newtheorem{theorem}{Theorem}[section]

\theoremstyle{definition}
\newtheorem{definition}[theorem]{Definition}

\newtheorem{proposition}[theorem]{Proposition}

\theoremstyle{remark}
\newtheorem{remark}[theorem]{Remark}

\numberwithin{equation}{section}

\theoremstyle{remark}

\newcommand{\bZ}{\mathbb{Z}}

\definecolor{gr90}{gray}{0.90}

\definecolor{gr75}{gray}{0.75}

\definecolor{gyblue}{cmyk}{0,0.5,0,0}

\newsavebox{\myonesquare}
\savebox{\myonesquare}{\textcolor{gr75}{\rule{12.5pt}{12.5pt}}}
\newsavebox{\mytwosquare}
\savebox{\mytwosquare}{\textcolor{gr75}{\rule{25pt}{12.5pt}}}
\newsavebox{\mythreesquare}
\savebox{\mythreesquare}{\textcolor{gr75}{\rule{37.5pt}{12.5pt}}}
\newsavebox{\myfoursquare}
\savebox{\myfoursquare}{\textcolor{gr75}{\rule{50pt}{12.5pt}}}

\newsavebox{\othertwosquare}
\savebox{\othertwosquare}{\textcolor{gr75}{\rule{25pt}{12.5pt}}}

\newcommand{\sqone}{\usebox{\myonesquare}}
\newcommand{\sqtwo}{\usebox{\mytwosquare}}
\newcommand{\sqthree}{\usebox{\mythreesquare}}

\newcommand{\sqother}{\usebox{\othertwosquare}}

\newcommand\strongof[1]{{#1}^+}
\newcommand \partitionof[1]{\widetilde{#1}}
\newcommand\reverse[1]{{#1}^*}

\newcommand{\xx}{\mathbf{x}}

\newcommand{\qschur}{\mathcal{S}}
\newcommand{\dzatom}{\mathcal{A}}
\newcommand{\dzchar}{\kappa}

\newlength{\cellsize} 
\newsavebox{\cell}

\sbox{\cell}{\begin{picture}(18,18) \put(0,0){\line(1,0){18}}
\put(0,0){\line(0,1){18}} \put(18,0){\line(0,1){18}}
\put(0,18){\line(1,0){18}}
\end{picture}}

\newcommand\cellify[1]{\def\thearg{#1}\def\nothing{}%
\ifx\thearg\nothing \vrule width0pt height\cellsize depth0pt\else
\hbox to 0pt{\usebox{\cell} \hss}\fi%
\vbox to \cellsize{ \vss \hbox to \cellsize{\hss$#1$\hss} \vss}}

\newcommand\tableau[1]{\vtop{\let\\\cr
\baselineskip -16000pt \lineskiplimit 16000pt \lineskip 0pt
\ialign{&\cellify{##}\cr#1\crcr}}}

\newcommand\bas[1]{\omit \vbox to \cellsize{ \vss \hbox to \cellsize{\hss$#1$\hss} \vss}}

\newcommand{\emt}{\mbox{ }}

\begin{document}

% \title[short text for running head]{full title}
\title{Refinements of the Littlewood-Richardson Rule}

%    Only \author and \address are required; other information is
%    optional.  Remove any unused author tags.

%    author one information
% \author[short version for running head]{name for top of paper}
\author{J. Haglund}
\address{Department of Mathematics, University of Pennsylvania, Philadelphia, PA 19104-6395, USA}
%\curraddr{}
\email{jhaglund@math.upenn.edu}
\thanks{}

%    author two information
\author{K. Luoto}
\address{Department of Mathematics, University of British Columbia, Vancouver, BC V6T 1Z2, Canada}
%\curraddr{}
\email{kwluoto@math.ubc.ca}
\thanks{}

%    author three information
\author{S. Mason}
\address{Department of Mathematics, University of California at San Diego, San Diego, CA 92093, USA}
%\curraddr{}
\email{skmason@math.ucsd.edu}
\thanks{}

%    author four information
\author{S. van Willigenburg}
\address{Department of Mathematics, University of British Columbia, Vancouver, BC V6T 1Z2, Canada}
%\curraddr{}
\email{steph@math.ubc.ca}
\thanks{}
 
%    \subjclass is required.
\subjclass[2000]{Primary 05E05; Secondary 05E10, 33D52}
%    The 2010 edition of the Mathematics Subject Classification is
%    now available.  If you are citing a classification from the
%    new scheme, use the following input coding instead.
%\subjclass[2010]{Primary }

\date{}

\dedicatory{}

\begin{abstract}

In the prequel to this paper, we showed how results of 
Mason  involving a new combinatorial formula
for polynomials that are now known as Demazure atoms 
(characters of quotients of Demazure modules, called standard bases by 
Lascoux and Sch\"utzenberger) could be used to define a new basis for the
ring of quasisymmetric functions we call ``Quasisymmetric Schur functions" (QS functions for short).  
In this paper we develop the combinatorics of these polynomials futher, by showing that the 
product of a 
Schur function and a Demazure atom has a positive expansion in terms of Demazure atoms.  
We use these techniques, together with the fact that both a QS function and a Demazure
character have explicit expressions as a positive sum of atoms,
to obtain the 
expansion of a product of a Schur function with
a QS function (Demazure character) as a positive sum of QS functions (Demazure characters).
Our formula for the coefficients in the expansion of a product of a Demazure character
and a Schur function into Demazure characters 
is similar to known results and includes in particular the
famous Littlewood-Richardson rule for the expansion of a product of Schur functions in terms of the
Schur basis.

\end{abstract}

\maketitle

%\noindent\emph{keywords: {key polynomials, nonsymmetric Macdonald polynomials, Littlewood-Richardson rule, quasisymmetric functions, Schur functions, tableaux}}

\section{Introduction}

A composition (weak composition) with $n$ parts is a sequence of $n$ positive 
(nonnegative) integers, respectively.  A partition is a composition whose parts are monotone
nonincreasing.  If $\tau$ is a weak composition, composition, or partition, we let
$\ell (\tau)$ denote the number of parts of $\tau$.  Throughout this article $\gamma$ is a weak composition
with $\ell (\gamma)=n$ while $\beta$ and $\lambda$ denote compositions and partitions, respectively,
with $\ell (\beta) \le n$, $\ell (\lambda) \le n$.  The polynomials in this paper
(Schur functions, Demazure atoms and characters, QS functions) depend on a finite set of
variables $X_n=\{x_1,x_2,\ldots ,x_n\}$ which we often omit for the sake of readability.

Symmetric functions in a set of variables $X_n$ play a central role in representation theory, and in recent
years have found increasing utility in several other branches of mathematics and physics such as
special functions, algebraic geometry, and statistical mechanics.  One of the most general
symmetric functions is a family of orthogonal polynomials $J_{\mu}(X_n;q,t)$ introduced by
Macdonald \cite{Mac88}, \cite{Mac} in $1988$, which depend not only on $X_n$ but also on a 
partition
$\mu$ and two extra parameters $q,t$.  The $J_{\mu}$ contain many of the most useful symmetric
functions as limiting or special cases.  In $1995$ Macdonald \cite{Mac96} introduced a very general family of
orthogonal polynomials called the nonsymmetric Macdonald polynomials $E_{\gamma}(X_n;q,t)$ which, 
although not symmetric functions, 
satisfy versions of most of the nice analytic and algebraic properties of the $J_{\mu}$.  
Macdonald showed how to
express $J_{\mu}$ as a linear combination of the $E_{\gamma}$, which can thus be thought of as more 
fundamental building blocks. 
Macdonald's defintion of the $E_{\gamma}$ was rather indirect, but 
in \cite{HHL} a new combinatorial formula for the (type A)
$E_{\gamma}$ was introduced. 
By letting $q=t=0$ and $q=t=\infty$ in this formula
we obtain new combinatorial formulas for Demazure characters (first studied by 
Demazure in \cite{Dem}) and Demazure atoms (called standard bases
by Lascoux and Sch\"utzenberger \cite{LaSc}), respectively.  These formulas are described in terms of
{\it skyline fillings},
which are combinatorial objects related to tableaux.  
Mason \cite{mason-1},\cite{mason-2} showed that 
many of the interesting 
properties of Demazure characters and atoms can be explained via the combinatorics of
skyline fillings.  In particular she developed a refinement of the well-known Robinson-Schensted-Knuth algorithm,
involving skyline fillings and weak compositions, which shows bijectively that 
the Schur function 
$s_{\lambda}(X_n)$ is a sum of those atoms corresponding to weak compositions $\gamma$ with $n$ parts
whose nonzero parts are a rearrangement of the parts of $\lambda$.  

One natural question to ask is how this decomposition of $s_{\lambda}$ into atoms 
compares with the well-known fact \cite[p. 361]{ECII}
that $s_{\lambda}$ is a sum, over standard Young tableaux $T$ of shape $\lambda$,
of Gessel's fundamental quasisymmetric
function $F_{\text{des}(T)}$.
In \cite{HLMvW} the
authors showed that, if $\strongof \gamma$ is the composition obtained by removing all
zero parts from $\gamma$ (so for example, $\strongof {100203401} = 12341$) then the sum of Demazure atoms,
over all $\gamma$ with $\strongof \gamma$ equaling a fixed composition $\beta$, is a sum of
certain fundamental quasisymmetric functions, and hence also quasisymmetric.  We call this sum 
the quasisymmetric Schur function (QS for short), denoted $\qschur _{\beta}(X_n)$
and note that $s_{\lambda}(X_n)$ is the sum, over all
compositions $\beta$ whose parts are a rearrangement of the parts of $\lambda$ (denoted $\tilde{\beta}=\lambda$), of $\qschur _{\beta}(X_n)$.  
In general there are fewer terms 
in this expansion than the expansion into Gessel's $F$'s; 
for example, if $\lambda$ is a rectangle, then there is only
one multiset permutation of the parts of $\lambda$ and hence $s_{\lambda} = \qschur _{\lambda}$.

The family of QS functions forms a new basis for the ring of quasisymmetric functions.
Although the product of two fundamental quasisymmetric functions expands as a positive sum of
fundamental quasisymmetric functions \cite{Gessel}, it turns out that the product of two QS functions
does not expand as a positive sum of QS functions.  In \cite{HLMvW} the authors showed though that 
if you multiply a QS function by either
a complete homogeneous symmetric function or an elementary symmetric function the result is a positive
sum of QS functions, which can be thought of as a version of the famous Pieri rule.
The current investigation grew out of an observation of the authors 
that the product of a Schur function and a QS function is a positive sum 
of QS functions.  Efforts to understand the coefficients in this expansion combinatorially led to the
discovery that the product of a Schur function and a Demazure atom has a positive expansion into
atoms, and that the coefficients in this expansion can be described in terms of analogues of 
Littlewood-Richardson tableaux (also known as Yamanouchi tableaux), in the context of skyline fillings.  
We prove this in
Section  \ref{sec:LR-atoms}, borrowing many ideas contained in the proof in Fulton's book \cite{Fulton} of the 
classical Littlewood-Richardson 
rule, replacing statements about semi-standard Young tableaux (SSYT) 
by corresponding statements about skyline fillings.
In Sections \ref{LRQS} and \ref{sec:LR-chars} we show how our Littlewood-Richardson rule for atoms 
leads to corresponding rules for
both QS functions and Demazure characters.  Since Schur functions are special cases of Demazure characters,
we obtain the classical Littlewood-Richardson rule as a special case.  Note that in \cite{ReSh} 
Reiner and Shimizono obtain a number of
results involving the expansion of various generalizations of skew Schur functions as a 
positive sum of Demazure characters, which yield identities similar in spirit 
to our expansion of the product of
a Schur function and a Demazure character.

Every Schubert polynomial can be written as a positive sum of type A Demazure characters.  This means that our results provide a method for expanding the product of an arbitrary Schubert polynomial and Schur function (in the same set of variables) as a positive sum of Demazure characters.   Combinatorial descriptions of the coefficients that arise in the product of a Schur polynomial and certain Schubert polynomials when expanded as a sum of Schubert polynomials are given in \cite{Kogan} and \cite{Kohnert}. Their proofs involve concepts that are similar to ones we use, and one wonders whether it is possible to recover these results, or generalizations of them, using our techniques.  One also wonders whether there is an underlying structure unifying our positivity expansions, and whether this structure is related to other structures underlying algebraic positivity, for example, Polo's notion of $B$ modules with excellent filtration \cite{Polo}, \cite{vdk}.

%%%%%
\section{Basic definitions and notation}

\subsection{Skyline diagrams}

\begin{figure}
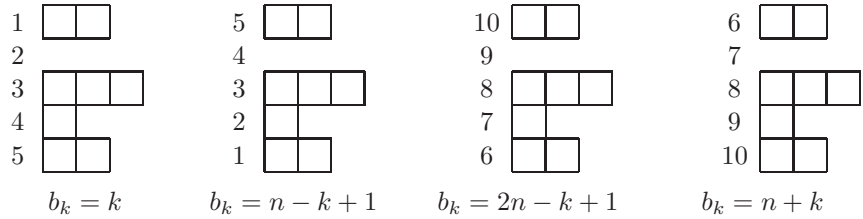

\begin{center}
\[
\begin{array}{cccc}
\tableau{
   \bas{1} & \emt & \emt \\
   \bas{2} \\ \bas{3} & \emt & \emt & \emt \\
   \bas{4} & \emt \\ \bas{5} & \emt & \emt 
}
\quad & \quad
\tableau{
   \bas{5} & \emt & \emt \\
   \bas{4} \\ \bas{3} & \emt & \emt & \emt \\
   \bas{2} & \emt \\ \bas{1} & \emt & \emt 
}
 \quad & \quad \quad
\tableau{
   \bas{10} & \emt & \emt \\
   \bas{9} \\ \bas{8} & \emt & \emt & \emt \\
   \bas{7} & \emt \\ \bas{6} & \emt & \emt 
}
 \quad & \quad \quad
\tableau{
   \bas{6} & \emt & \emt \\
   \bas{7} \\ \bas{8} & \emt & \emt & \emt \\
   \bas{9} & \emt \\ \bas{10} & \emt & \emt 
}
 \vspace{6pt} \\
 b_k=k  & b_k=n-k+1 & b_k=2n-k+1 & b_k=n+k\\
\end{array}
\]
\caption{Skyline diagrams with $n=5$, composition $(2,0,3,1,2)$, and four typical basements}
\label{fig-skyline}
\end{center}
\end{figure}

A \emph{skyline diagram} 
is a collection of boxes, or cells, arranged into left-justified rows  
\footnote{This differs 
slightly from
the convention in \cite{HHL}, \cite{mason-1}, \cite{mason-2}, where skyline diagrams are
arranged in bottom-justified columns.}.
To each skyline diagram we associate a weak composition,
whose $k$th part is the number of cells in the $k$th row of 
the diagram, where the top row is viewed as row $1$, the row below it row $2$, et cetera.
Skyline diagrams are augmented by a \emph{basement}, an extra column on the 
left (considered to be the $0$-th column) containing positive integers.  We let $b_k$ denote the
entry in the
$k$th row of the basement.    
In most of our examples the basement will either satisfy $b_k=k$, $b_k=n-k+1$, $b_k=n+k$, 
or $b_k=2n-k+1$ for $1\le k \le n$,
as in the diagrams in Figure \ref{fig-skyline}. 

Let $\gamma,\delta$ be weak compositions with $\gamma \subseteq \delta$, 
i.e. $\gamma_i \le \delta_i$  for $1\le i \le n$.

A \emph{skew skyline diagram} of shape $\delta/\gamma$ is obtained by starting with a skyline diagram of shape $\delta$ with basement values $(b_1,b_2,\ldots,b_n)$, and considering the cells of $\gamma$ to be an extension of the basement, placing the value $b_k$ in each of the cells in the $k$-th row of $\gamma$.
A skyline diagram of shape $\delta$ can naturally be 
viewed as a skew skyline
diagram of shape $\delta/ (0,0,\ldots,0)$.  Note that if the parts of $\delta$ and of $\gamma$ are
monotone decreasing, and we remove the basement, we get a \emph {skew} Ferrers shape.

A skyline filling (skew skyline filling) is an assignment of positive integers to the
cells of a skyline (skew skyline) diagram, respectively.
Central to our constructs involving skyline fillings is a \emph{triple} of cells, 
of which there are two types.
A \emph{type A triple} in a diagram of shape $\delta /\gamma$ is a set of three cells $a,b,c$ of the 
form $(i,k),(j,k),(i,k-1)$ for some pair of rows $i < j$ of the diagram and some 
column $k>0$, where row $i$ is at least as long as row $j$, i.e. $\delta_i \geq \delta_j$.
A \emph{type B triple}  is a set of three cells $a,b,c$ of the form $(j,k+1),(i,k),(j,k)$ for 
some pair of rows $i < j$ of the diagram and some column $k\geq 0$, where row $i$ is 
strictly shorter than row $j$, i.e. $\delta_i < \delta_j$.
Note that basement cells can be elements of triples.
As noted below, in this article our
fillings have weakly decreasing row entries left-to-right, so we always have the entry values $c \geq a$.
We say that a triple of either type is an \emph{inversion} triple if the relative order of the entries is
either $b < a \leq c$ or $a \leq c < b$.
Otherwise we say that the triple is a \emph{coinversion} triple, i.e. $a \leq b \leq c$.

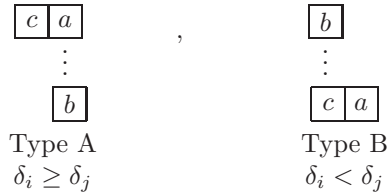
\[
\begin{array}{ccc} \vspace{6pt}
\tableau{   c & a   \\  & \bas{{\vdots}} } 
&  \qquad , \qquad  &
\tableau{   b   \\  \bas{{\vdots}}} 
 \\
 \tableau{ & b} & & \hspace*{.2in} \tableau{c & a} \\
\mbox{Type A} & & \hspace*{.2in} \mbox{Type B} \\
\delta_i \geq \delta_j & & \hspace*{.2in} \delta_i < \delta_j
\end{array}
\]

\medskip
\noindent
A {\it semistandard skyline filling} (SSK) is a (skew) skyline filling where
\begin{enumerate}[(i)]
\item 
each row is weakly decreasing left-to-right (including the basement), and 
\item  all triples (including triples with cells in the basement) are inversion triples.
\end{enumerate}

\begin{remark}
\label{large}
Note that since basement values are constant across rows, for any choice of basement values any triple 
involving three basement cells is forced to be an inversion triple.
Furthermore, if we have a skew skyline diagram with basement $b_k=2n-k+1$, all entries in the
basement are larger than $n$, the biggest entry outside the basement.
Therefore, the actual
values of the $b_k$ are not relevant, as long as they are decreasing from top to bottom and all larger than $n$.
For this reason we often draw the basement $b_k=2n-k+1$ with ``$*$" symbols in place of the $b_k$, where we
think of the $*$ as a virtual $\infty$ symbol, larger than any entry,  and we refer to this basement
as the \emph{large basement}.  To determine whether a triple involving two $*$ symbols is an inversion triple
or not, we view $*$ symbols in the same row as being equal, and $*$ symbols in a given column as
decreasing from top to bottom.  In our identities involving the large basement and polynomials depending 
on $X_n$, we can let $n \to \infty$ and view the identity as holding in the infinite set of
variables $X = \{x_1,x_2,\ldots\}$.
%One other thing to note is that columns of height
%zero in our skew skyline shape which occur after all the columns of nonzero height have no effect on 
%whether a filling is an SSK or not, so when drawing our diagrams we can leave 
%such columns out.
\end{remark}

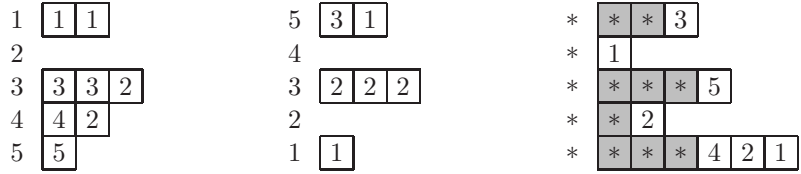
\begin{figure}
\begin{center}
\[
\begin{picture}(400,72)
\put(0,72){\tableau{
   \bas{1} & 1 & 1 \\
   \bas{2} \\ 
   \bas{3} & 3 & 3 & 2 \\
   \bas{4} & 4 & 2 \\ 
   \bas{5} & 5 
}}
\put(150,72){\tableau{
   \bas{5} & 3 & 1 \\
   \bas{4} \\ 
   \bas{3} & 2 & 2 & 2 \\
   \bas{2}  \\ 
   \bas{1} & 1 
}}
\put(323,72){\sqtwo}
\put(323,36){\sqthree}
\put(323,18){\sqone}
\put(323,0){\sqthree}
\put(300,72){\tableau{
   \bas{*} & * & * & 3 \\
   \bas{{*}} & 1 \\ 
   \bas{*} & * & * & * & 5 \\
   \bas{{*}} & *  & 2 \\ 
   \bas{{*}} & * & * & * & 4 & 2 & 1
}}
\end{picture}
\]
\caption{SSK of shapes $(2,0,3,2,1)$, $(2,0,3,0,1)$, 
and $(3,1,4,2,6)/(2,0,3,1,3)$}
\label{fig-SSK}
\end{center}
\end{figure}

Figure \ref{fig-SSK} gives examples of SSK for various shapes $\delta/\gamma$.  
The portion of the
basement whose cells are part of $\gamma$ are indicated by light gray backgrounds.
It is shown in \cite{mason-1} that every $\text{SSK}$ is \emph{non-attacking}, 
meaning that the entries within each column are all distinct, 
and that two cells $a=(i,k)$ and $b=(j,k+1)$ can only have the same value if $i\geq j$.  

%\medskip \noindent
%A \emph{semistandard composition tableau} (CT) is a composition filling where
%\begin{enumerate}[(i)]
%\item the first column is strictly increasing, top-to-bottom,
%\item each row is weakly decreasing left-to-right, and
%\item  all triples are inversion triples.
%\end{enumerate}

\setlength{\fboxrule}{1pt}%
\setlength{\fboxsep}{3pt}%

\subsection{Contretableaux and reading words}

A \emph{contretableau} (CT) is a Ferrers shape filled with positive integers
where the entries within each row \emph{decrease} weakly left-to-right and 
the entries within each column \emph{decrease} strictly top-to-bottom.  
We let $\text{CT}(n)$ denote the
set of CT with entries from the set $[n]=\{1,2,\ldots ,n\}$.  Note $\text{CT}(n)$ is trivially
in bijection with $\text{SSYT}(n)$, the set of $\text{SSYT}$ with entries from $[n]$,
by applying the map $j \to n-j+1$ to each entry of a given CT.
%A \emph{standard} contretableau (SCT) is a CT in which the set of entries is exacty $[n]$.  

Since CT are trivially equivalent to SSYT, it is no surprise that
all of the concepts, definitions, operations (such as insertion and evacuation), 
propositions, and theorems regarding SSYT have CT-counterparts, 
and the proofs of such results are completely analogous.  We include in this 
section several of the classical notions most pertinent to our results; 
the (SSYT versions of the) fully developed theory can be found in \cite{Fulton} or \cite{ECII}.

The \emph{row reading order} of a (possibly skew) skyline diagram or Ferrers shape
is a total ordering of the cells 
where $(i,j) <_{row} (i',j')$ if either $i > i'$ or ($i=i'$ and $j < j'$).
That is, the row reading order reads the cells  left-to-right in each row, 
starting with the bottommost row and proceeding upwards to the top row, ignoring basement entries if they exist.
The \emph{row word} of a filling $T$, denoted $\text{rowword}(T)$ is the sequence of integers 
formed by the entries of $T$ taken in row reading order.

We also use a slightly different reading order on diagrams, which we refer to as the \emph{column reading order}.
In the column reading order, we have $(i,j) <_{col} (i',j')$ if either $j > j'$ or ($j = j'$ and $i < i'$).
That is, the column reading order reads the cells from top to bottom within each column, 
starting with the rightmost column and working leftwards, again ignoring any basement entries.
The \emph{column word} of a filling $T$, denoted $\text{colword}(T)$, is the sequence of 
integers formed by the entries of $T$ taken in column reading order.  For example, for the rightmost
SSK in Figure \ref{fig-SSK}, the row word is $4212513$ and the column word  
is $1254321$.

%%%%%

\begin{figure}[hb]
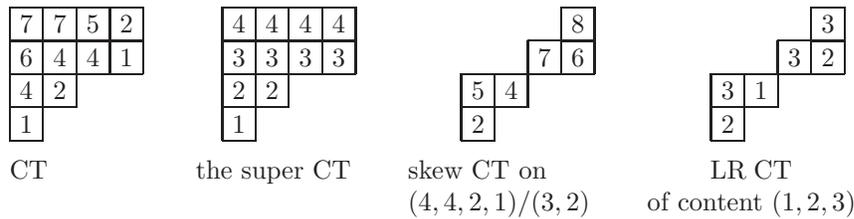

\begin{center}
\[
\begin{array}{lllc}
\tableau{
    7 &   7 & 5 & 2  \\
    6 &  4  & 4 & 1  \\
    4 &  2   \\
    1   \\
} \quad & \quad
\tableau{
    4 &   4 & 4 & 4  \\
    3 &  3 & 3 & 3  \\
    2 &  2   \\
    1   \\
} \quad & \quad \quad
\tableau{
    &  & & 8  \\
    &   & 7 & 6  \\
    5 &  4   \\
    2   \\
} \quad & \quad \quad
\tableau{
    &  & & 3  \\
    &   & 3 & 2  \\
    3 &  1   \\
    2   \\
} \vspace{6pt} \\
 \text{CT} & \text{the super CT} & \text{skew CT on} & \text{LR} \;\text{CT} \\
  & & (4,4,2,1)/(3,2) & \text{of content $(1,2,3)$}
\end{array}
\]
\caption{CT examples}
\label{fig-contretableau}
\end{center}
\end{figure}

\begin{definition}
For a word $w=w_1w_2\cdots w_n$ (or sequence $(w_1,\ldots, w_n)$) we let 
$w^{*}$ denote the reverse word $w_nw_{n-1}\cdots w_2w_1$ (or sequence $(w_{n},\ldots ,w_1)$).
The \emph {content} of $w$ is the sequence $(c_1,c_2,\ldots ,c_n)$ where $c_i$ is the number of occurrences
of $i$ in $w$.  The content of a CT $T$ is the content of $\text{colword}(T)$.
\end{definition}

For a partition $\lambda$ with $r$ parts,
the \emph{super CT} of shape $\lambda$ (denoted $U_\lambda$) is the unique CT with content 
$\lambda^{*}$, 
as in Figure \ref{fig-contretableau}.
For the``contre-'' analog $V\gets W$ of inserting/recording a biword $W$ into a 
CT $V$, the biletters of $W$ are sorted into \emph{reverse} (i.e. weakly decreasing) lexicographical order.

A lattice word is a word (or sequence) $w=w_1w_2\cdots w_n$
where in any initial segment
$w_1w_2\cdots w_i$ there are at least as many occurrences of the number $j$ as $j+1$, for each $j\ge 1$.
We say $w$ is \emph {contre-lattice} if
in any initial segment there are at least as many occurrences of the number $j$ as $j-1$, for each $1< j \leq r$, where
$r$ is the maximum of the $w_i$.
We say a word or sequence $w$ is \emph {regular contre-lattice}
if it is contre-lattice  
and the minimum of the $w_i$ is $1$.
We define a \emph {Littlewood Richardson skew CT} to be a skew CT the reverse of whose row reading word is regular 
contre-lattice.
We often abbreviate ``Littlewood-Richardson" by ``LR".

\begin{proposition} \cite[Section $5.2$]{Fulton}
\label{prop:LRcontretabs}
Let $S$ be a skew CT with content $\mu^{*}$.
Then the following are equivalent.
\begin{enumerate}[(i)]
\item $S$ is an LR skew CT, i.e. $\text{rowword}(S)^{*}$ is 
a regular contre-lattice word.
\item $\text{colword}(S)$ is a regular contre-lattice word.
\item $\text{rect}(S) = U_\mu$, the super CT of shape $\mu$.  
(Here $\text{rect}(S)$ is the ``rectification" of $S$ - see \cite{Fulton}.)
\end{enumerate}
\end{proposition}

\subsubsection{Combinatorial formulas}

Recall the well-known combinatorial formula for the Schur function 
\begin{equation}  \label{eq:schurCT}
s_\lambda = \sum_{T\in \text{SSYT}(n),\atop \text{shape}(T)=\lambda} \xx^T.
\end{equation}
The following combinatorial formulas for Demazure atoms $\dzatom _{\gamma}$ and Demazure 
characters $\dzchar _{\gamma}$ follow as limiting cases of results in \cite{HHL} 
\begin{eqnarray}
 \label{eq:atom} \dzatom_\gamma &= &  \sum_{Y\in \text{SSKI}(n),\atop \text{shape}(Y)=\gamma} \xx^Y  \\
 \dzchar_\gamma &= &  \sum_{Y\in \text{SSKD}(n),\atop \text{shape}(Y)=\gamma ^{*}} \xx^Y 
\end{eqnarray}
where $\text{SSKI}(n)$ is the set of all SSK with basement $b_k=k$ and entries in $[n]$,
and $\text{SSKD}(n)$ is the set of all SSK with entries in $[n]$ and $b_k=n-k+1$.

Taylor~\cite{Taylor}, building on the work of \cite{ABW}, \cite{LM}, \cite{ReSh-1}, \cite{ReSh-2}, \cite{ReSh-3}, describes a class of generalized tableaux over signed alphabets which he calls straight. Straight tableaux have overall partition shape and include the classical skew SSYT for all-positive alphabets. They are defined in part using a triple condition on cells which, if we reverse his inequalities and consider only positive alphabets, would be the same as our type A inversion triple condition. Additionally, those SSKD with partition shape are equivalent to CT and thus to SSYT, which are special cases of straight tableaux. It would be interesting to know whether there is a more significant relationship between Taylor's straight
tableaux and the objects in this article.

\subsubsection{A bijection between $\text{SSKI}(n)$ and $\text{CT}(n)$}

There exists a simple bijection $\rho$ between $\text{SSKI}(n)$ and $\text{CT}(n)$~\cite{mason-1}. 
Given $Y\in \text{SSKI}(n)$, 
one obtains the corresponding 
CT by sorting the entries within each column, as in the example below.

\[
\begin{array}{ccc}
\tableau{
   \bas{1} & 1 & 1 \\
   \bas{2} \\ \bas{3} & 3 & 3 & 2 \\
   \bas{4} & 4 & 2 \\ \bas{5} & 5 
}
\quad & \xrightarrow {\rho} \quad
\tableau{ 5 & 3 & 2 \\ 4 & 2\\ 3 & 1\\ 1}
\end{array}
\]

The inverse $\rho^{-1}$ is only slightly more intricate.  Given $T\in \text{CT}(n)$, 
map the leftmost column of $T$ into an empty element of $\text{SSKI}(n)$ by placing the entries of this column, 
beginning with the largest, into the highest row of the leftmost column whose rightmost entry is weakly greater.  
Repeat this procedure with each of the remaining columns.  
% Note that the map $\rho$ and its inverse preserve the set of column entries. 
One important property to note about the bijection is that it preserves the set of entries 
within each column.  (We say set as opposed to multiset since all of our tableau-like structures 
require that all entries within a column be distinct.)

\subsubsection{Pieri rules}  \label{sec:Pieri}
Pieri rules for multiplying a QS function by a complete homogeneous symmetric function $s_k$ 
or an elementary symmetric function $s_{1^k}$ 
are presented in \cite{HLMvW}.  
By the same method one can derive Pieri rules for multiplying a Demazure atom by either $s_k$ or $s_{1^k}$.
The intersecting case $s_1$, the ``single box case'', can be described as follows.
Given a weak composition $\delta$ containing a part with value $k$, $k> 0$, define $rem_k(\delta)$ to be the weak composition resulting from decrementing the last (rightmost) part of $\delta$ that has value $k$.
For example, 
\[  rem_2(1,0,4,2,0,1,2,3) = (1,0,4,2,0,1,1,3). \]
We likewise define $rem_k(\beta)$ for compositions, where the result is collapsed to a 
composition by removing any resulting zero part.
For example, \[  rem_1(1,4,2,1,2,3) = (1,4,2,2,3). \]
Now the ``single box'' Pieri rule can be described as
\begin{eqnarray}
\dzatom_\gamma (X_n) \cdot s_1(X_n) &=& \sum_\delta \dzatom_\delta (X_n)\\
\qschur_\alpha (X_n) \cdot s_1 (X_n) &=& \sum_\beta \qschur_\beta (X_n) 
\end{eqnarray}
where $\delta$ runs over all weak compositions satisfying 
$\gamma = rem_k(\delta)$ for some positive integer $k$, 
and similarly $\beta$ runs over all compositions satisfying $\alpha = rem_k(\beta)$ for 
some positive integer $k$.  Note the close similarity between these rules and the corresponding rule for
Schur functions \cite[p. 337]{ECII}.

%%%%%%
\section{Properties of skyline fillings}  \label{sec:sky-props}

%We prove a some intermediary propositions regarding a mild generalization of LRS's.
For a given cell $x$ in a skyline diagram, we let $\text{row}(x)$ denote the row containing $x$.
Say that an SSK $Y$ on any basement 
is \emph{contre-lattice} if its column reading word is contre-lattice.
%We note that for any skyline the elements within any given column are all distinct.
Suppose $Y$ is an SSK on any basement, where $Y$ has $t$ columns.
For all $1 \leq k \leq t$, let $C_k$ be the set of entries in column $k$ of $Y$, 
excluding basement entries.
Call these the \emph{column sets} of $Y$.
Sort each $C_k$ into decreasing order and form the word $w_Y = C_t C_{t-1}\cdots C_2 C_1$.
We say that $Y$ is \emph{loosely contre-lattice} if $w_Y$ is contre-lattice.

%%%
\begin{proposition} \label{prop:r-contre}
Let $Y$ be an element of $\text{SSK}(n)$ on any basement satisfying $b_k>n$ for $1\le k \le n$.
Then $Y$ is contre-lattice if and only if $Y$ is loosely contre-lattice.
\end{proposition}
\begin{proof}
Assume $Y$ is contre-lattice and $Y$ has $t$ columns.
Let $C'_k$ be the sequence of the elements of the $k$th column of $Y$ in column reading order, 
so that $\text{colword}(Y) = C'_tC'_{t-1}\cdots C'_2C'_1$. 
By assumption, $\text{colword}(Y)$ is contre-lattice.
If within this word we transpose any adjacent pair $w_i w_j$ of letters in the word, 
where $w_i < w_j$, then the resulting word retains the contre-lattice property.
In particular, if we sort each of the $C'_k$ into decreasing order to obtain $C_k$, 
the resulting word $w_Y = C_t C_{t-1}\cdots C_2 C_1$ retains the contre-lattice property.
Thus $Y$ is loosely contre-lattice.

Conversely, assume that $Y$ is loosely contre-lattice.
Label the cells of $Y$ according to their contents and place in the column reading order of $Y$.
Specifically, we identify a cell of $Y$ as $x_j$ when the cell contains the $j$th 
occurrence of the entry $x$ in $\text{colword}(Y)$.
Let $m$ be the number of $r$'s in $\text{colword}(Y)$, and for each $1\leq k \leq m$ 
let $S_k = \{x_k : 1\leq x \leq r, \, x_k \in Y\}$.
To show that $Y$ is contre-lattice,
it suffices to show that for each $k$, the cells of $S_k$, as they appear in the 
column reading order of $Y$, are in strictly decreasing order of their contents.
Since $Y$ is loosely contre-lattice, if two cells of $S_k$ are in different columns of $Y$, 
then they appear in the column reading order of $Y$ in strictly decreasing order of their contents.
Thus it suffices to show that if two cells of $S_k$ lie in the same column of $Y$ then they also 
appear in strictly decreasing order of their contents.
That is, we need to show that if $x < y$ and $x_k,y_k \in S_k$ are in the same column, 
then $y_k$ appears above $x_k$ in that column.

Seeking a contradiction, suppose this is not the case.
Among all such violating pairs of values, choose $x$ and $y$ such that $|y-x|$ is minimized.

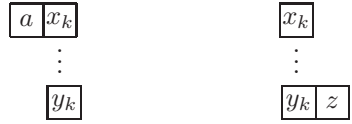
\[
\begin{array}{ccc} \vspace{6pt}
\tableau{   a & x_k   \\  & \bas{{\vdots}} } 
&  \qquad  \qquad  &
\tableau{   x_k  \\  \bas{{\vdots}}} 
 \\
 \tableau{ & y_k} & & \hspace*{.2in} \tableau{y_k & z}
\end{array}
\]

We consider two cases.
In the first case, $\text{row}(x_k)$ is at least as long as $\text{row}(y_k)$, corresponding to the 
first diagram above,
where $a$ is the entry immediately to the left of $x_k$.
Without loss of generality, we may assume that $k$ is largest among such indices for this case.
Since $Y$ has no coinversion triples, it must be the case that $x_k \leq a < y_k$,
implying that $a$ is not in the basement.
If $a = x$, then the cell with entry $a$ shown is in fact $x_{k+1}$.
This implies that $y_{k+1}$ is also in the same column as $x_{k+1}$.
Since $Y$ is non-attacking, $y_{k+1}$ must appear weakly below $y_k$, and hence below $x_{k+1}$,
contrary to the assumption that $k$ is maximal.
Thus $x_k < a < y_k$. 
This in turn implies that $a_k$  is in the same column as $x_k$ and $y_k$ and also that $a = a_{k+1}$.
Since $Y$ is non-attacking, $a_k$ must appear weakly above $a_{k+1}$, and hence above $y_k$.
But then $y$ and $a$ form a violating pair with $|y-a|<|y-x|$, contrary to our 
assumption that $|y-x|$ is minimal.  Thus we have a contradiction in this case.

The other case is that $\text{row}(x_k)$ is shorter than $\text{row}(y_k)$, corresponding to the second diagram above, where $z$ is the entry immediately to the right of $y_k$.
Without loss of generality, we may assume that $k$ is smallest among such indices for this case.
Since $Y$ has no coinversion triples, it must be the case that $x_k < z\leq y_k$.
If $z = y$, then the cell with entry $z$ shown is in fact $y_{k-1}$.
This implies that $x_{k-1}$ is also in the same column as $y_{k-1}$.
Since $Y$ is non-attacking, $x_{k-1}$ must appear weakly above $x_k$, and hence above $y_{k-1}$,
contrary to the assumption that $k$ is minimal.
Thus $x_k < z < y_k$. 
This in turn implies that $z_k$  is in the same column as $x_k$ and $y_k$ and also that $z = z_{k-1}$.
Since $Y$ is non-attacking, $z_k$ must appear below $z_{k-1}$, and hence below $x_k$.
But then $x$ and $z$ form a violating pair with $|z-x|<|y-x|$, contrary to our assumption that $|y-x|$ is minimal.  Thus we have a contradiction in this case.
Thus in all cases we obtain a contradiction, which completes the proof.
\end{proof}

%%%
\begin{proposition} \label{prop:r-LRS-reduction}
Let $Y$ be a contre-lattice element of $\text{SSK}(n)$
on a decreasing (i.e. $b_i > b_{i+1}$ for
$1 \le i <n$) basement.
Let $x$ be the smallest entry value in $Y$, and let $x_1$ be the cell containing the 
rightmost entry of value $x$, i.e. the first $x$ in column reading order.  Then
 $x_1$ is the rightmost cell of its row, say row $i$, and for every row $i' >i$, row $i$ 
and row $i'$ have different lengths.
\end{proposition}
\begin{proof}
That $x_1$ is at the end of its row is immediate since $x$ is the smallest entry value in $Y$.
Suppose that there is some row $i' > i$ of the same length as row $i$.
Let $z$ be the entry in the last cell of row $i'$.
\[ \tableau{ \bas{b} & \bas{\cdots} & u & v & \bas{\; \, \cdots} & x_1 \\ \\ \bas{b'} & \bas{\cdots} & w & y & \bas{\; \, \cdots} & z } \]
Since $x$ is the smallest entry value in $Y$, we have $x < z$.
On the other hand, since the basement is decreasing, the basement entries for the rows are related by $b >b'$, where $b = b_i= Y(i,0)$ and $b'=b_{i'}=Y(i',0)$. 
Thus there must exist some column $j$ such that 
$v < y$ and $u >w$,  where $v = Y(i,j+1)$, $y = Y(i',j+1)$, $u = Y(i,j)$, and $w = Y(i',j)$.
This implies that $v<y<u$, which would form a type A coinversion triple, a contradiction.
Thus there can be no such row $i'$.
\end{proof}

%%%
\begin{proposition} \label{prop:r-any-reduction}
Let $Y$ be a contre-lattice element of $\text{SSK}(n)$ on any basement.
Let $x$ be the smallest entry value in $Y$, and let $x_1$ be the cell containing the 
rightmost entry of value $x$, i.e. the first $x$ in column reading order.  Then
the skyline diagram filling $Y' =Y-x_1$ obtained from $Y$ by simply removing 
cell $x_1$ is also a
contre-lattice SSK . 
\end{proposition}
\begin{proof}
As above, assume $x_1$ is in row $i$.  
Since $Y$ is already an $\text{SSK}$, to show that $Y'$  is an $\text{SSK}$,  
it suffices to show that removing $x_1$ does not introduce any coinversion triples, 
which could only happen between row $i$ and some other row $i'$.
No type B coinversion triples could be introduced since by Proposition \ref{prop:r-LRS-reduction}
 there are no rows in $Y$ 
below row $i$ of the same length as row $i$.
Any type A coinversion triples introduced would have to be between a row $i' < i$ of 
length one less than that of row $i$ in $Y$.  
Suppose that such a  conversion triple $u,v,w$ exists in $Y'$ between rows $i'$ and $i$, as shown.
\[ \tableau{ \bas{\cdots} & w & u & \bas{\; \, \cdots} & a & c & \bas{\; \, \cdots} & y 
\\ \\ \bas{\cdots} & &  v& \bas{\; \, \cdots} &  b & d & \bas{\; \, \cdots} & z & x_1} \]
The relation between these values must be $u < v < w$.
In particular, $u < v$.
On the other hand, the triple $x_1,y,z$ occurring at the end of rows $i$ and $i'$ in $Y$, 
as shown ($z$ is to the immediate left of $x_1$), must be an inversion triple, 
and since $x$ is the smallest entry value in $Y$, this implies the order $x \leq z < y$.
In particular, $y > z$.
Thus there must exist some column $j$ such that 
$a < b$ and $c >d$,  where $a = Y(i^\prime,j)$, $b = Y(i,j)$, $c = Y(i^\prime,j+1)$, and $d = Y(i,j+1)$.
This implies that $d < c \leq a<b$.
In particular  $d,a,b$ would form a type B coinversion triple in $Y$, a contradiction.
Thus there can be no such type A coinversion triple in $Y'$, and so $Y^{\prime}$ is an $\text{SSK}$.

Lastly, removing the first occurring smallest-value letter from a contre-lattice 
word clearly leaves another contre-lattice word, and so $Y'$ is also contre-lattice.
\end{proof}

\begin{proposition}{\label{!shape}}
Let $Y$  be a contre-lattice element of $\text{SSK}(n)$
of shape $\delta/\gamma$ on any basement of shape $\gamma$.
Let $\sigma$ be any permutation of $\delta$.
Then there exists a unique contre-lattice SSK $T$ of overall shape $\sigma$ on a large basement  $b_i=2n-i+1$ having the same column sets as $Y$.
Moreover,  $T$ has shape $\sigma/\tau$ for some basement shape $\tau$ that is a permutation of $\gamma$.
\end{proposition}
\begin{proof}
We proceed by induction on the number of cells in $Y$, that is, $|\delta/\gamma|$.
Start with the unfilled skyline diagram of shape $\sigma$.
Let $x$ be the smallest entry value in $Y$, and let $j$ be the column of the rightmost occurrence of $x$ in $Y$.  
Since $x$ is the smallest entry value, it occurs at the end of some row of $Y$ of length $j$.
Proposition \ref{prop:r-LRS-reduction} tells us that if $T$ exists, then the entry $x$ in column $j$ of $T$ must occur in the last row of $T$ of length $j$, say row $i$, which exists since $\sigma$ is a permutation of $\delta$.
Let $Y'$ be the SSK obtained by removing the cell in column $j$ of $Y$ that contains $x$.
Let $\sigma'$ be the shape obtained by removing the last cell of the last row of length $j$ in $\sigma$.
If $|\delta/\gamma|=1$, we are done; we set the basement of $T$ to be of shape $\sigma'$.

Otherwise, by Proposition \ref{prop:r-any-reduction}, $Y'$ is also a contre-lattice SSK, say of shape $\delta'/\gamma$, and clearly $\sigma'$ is a permutation of $\delta'$.
By our induction hypothesis there is a unique contre-lattice SSK $T'$ of overall shape $\sigma'/\tau$ on a large basement  $b_i=2n-i+1$ having the same column sets as $Y'$, where $\tau$ is a permutation of $\gamma$.  We want to show that if we append a cell containing $x$ to row $i$ of $T'$, which must be in column $j$, then the resulting filling $T$ is an SSK.
Since $x$ was the minimum entry of $Y$ and column $j$ its rightmost appearance, the cell added to $T'$ 
to form $T$ is the rightmost minimum entry of $T$. 
In particular, it is less than or equal to the entry to its immediate left, so row $i$ of $T$ is weakly decreasing, as are all other rows of $T$.

It remains to check the triple conditions.
Consider row $i'$, where $i'\neq i$.
If the relative order of the lengths of the two rows $i$ and $i'$ is unchanged when comparing $T'$ to $T$, then the type of triples between the two rows remains the same, and we only need consider any new triple formed by adding the new cell.
In any new triple formed, $x$ lies in row $i$ while the other two cells of the triple lie in row $i'$,
and since $x$ is the rightmost occurrence of the minimum value in $T$, it cannot form a coinversion triple between the two rows.

\[
\begin{array}{ccc}
\tableau{
    & & & \bas{j} \\
   \bas{i:} & &  \bas{\cdots} & x & & }
   
\quad \quad \quad & {\rm or} \quad \quad \quad
\tableau{
    & & & & \bas{j} \\
   \bas{i':} & & \bas{\cdots} & c & a & \bas{\; \, \cdots}}
\\ \\
\tableau{\bas{i':} & & \bas{\cdots} & c & a & \bas{\; \, \cdots} & & 
} \hspace*{.1in} & \hspace*{.2in} \tableau{ & \bas{i:} & & \bas{\cdots} & b & x }
\end{array}
\]

The only remaining case  is when $i'<i$ and $\sigma'_{i'} = \sigma'_i$,  when $\sigma_{i'} < \sigma_i$.
In this case, whereas $T'$ had type A triples between the two rows, now $T$ has type B triples between them.
Suppose that one of these type B triples in $T$ is a coinversion triple, say $v,w,u$ as in the diagram below, where $v$ and $w$ are in column $j'$, and where possibly the cell $u$ is the cell at the end of row $i$.

$$\tableau{& & & & & & & & \bas{j'} & & & & \bas{j}} $$
$$\tableau{\bas{i':} & \bas{b'} & \bas{\cdots} & a & c & \bas{\; \, \cdots} & v & & \bas{\cdots} & y} $$

$$\tableau{ & \bas{i:} & \bas{b} & \bas{\cdots} & b & d & \bas{\; \, \cdots} & w & u & \bas{\; \, \cdots} & z & x}$$

This requires that $u\leq v < w$.  
On the other hand, since $T$ and $T'$ share a common decreasing basement, the basement entries in these rows satisfy $b' > b$.
This implies that there exists some pair of adjacent columns in the range 0 to $j'$ inclusive containing the cells $a$, $b$, $c$, and $d$ of the two rows as shown such that $a>b$ and $c<d$.
But that would imply that $c<d<a$, forming a type A coinversion triple in $T'$, contrary to the fact that $T'$ is a valid SSK.
Thus all the type B triples between the two rows in $T$ are inversion triples.
In all cases,  $T$ is a valid SSK.

By Proposition \ref{prop:r-contre}, $Y$ is loosely contre-lattice.
Since $T$ has the same column sets as $Y$, $T$ is therefore also  loosely contre-lattice,
and again by Proposition \ref{prop:r-contre}, $T$ is contre-lattice.
\end{proof}

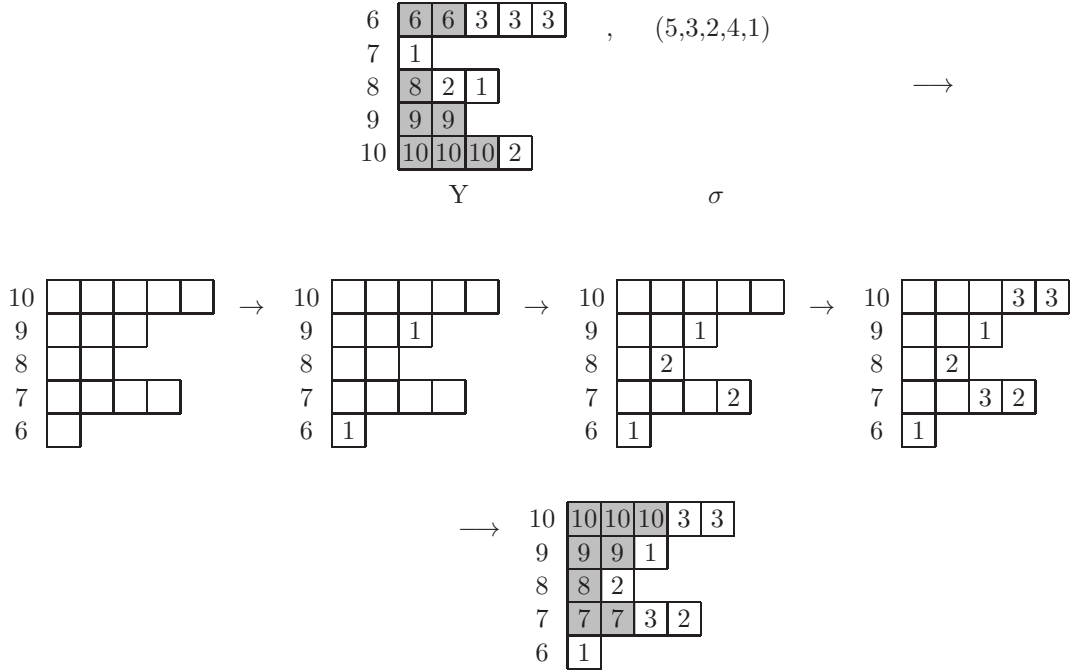
\begin{figure}[ht]
\begin{center}
\[
\begin{picture}(500,330)
\put(173,330){\sqtwo}
\put(173,294){\sqone}
\put(173,276){\sqtwo}
\put(173,258){\sqthree}
\put(150,330){\tableau{ \bas{{6}} & {6} & {6} & 3 & 3 & 3 \\
\bas{{7}} & 1 \\
\bas{{8}} & {8} & 2 & 1 \\
\bas{{9}} & {9} & {9} \\
\bas{{10}} & {10} & {10} & {10} & 2} \; \; , \; \;  (5,3,2,4,1)}
\put(450,300){$\longrightarrow$}
\put(200,240){Y}
\put(340,240){$\sigma$}
\put(-40,180){\tableau{
   \bas{10} & \emt & \emt & \emt & \emt & \emt \\
   \bas{9} & \emt & \emt & \emt \\ 
   \bas{8} & \emt & \emt   \\
   \bas{7} & \emt & \emt & \emt & \emt \\ 
   \bas{6} & \emt 
} \; $\rightarrow$ \;
\tableau{
   \bas{10} & \emt & \emt & \emt & \emt & \emt \\
   \bas{9} & \emt & \emt & 1 \\ 
   \bas{8} & \emt & \emt   \\
   \bas{7} & \emt & \emt & \emt & \emt \\ 
   \bas{6} & 1 
}  \; $\rightarrow$ \;
 \tableau{
   \bas{10} & \emt & \emt & \emt & \emt & \emt \\
  \bas{9} & \emt & \emt & 1 \\ 
   \bas{8} & \emt & 2   \\
   \bas{7} & \emt & \emt & \emt & 2 \\ 
   \bas{6} & 1 
} \; $\rightarrow$ \;
 \tableau{
   \bas{10} & \emt & \emt & \emt & 3 & 3 \\
   \bas{9} & \emt & \emt & 1 \\ 
   \bas{8} & \emt & 2   \\
   \bas{7} & \emt & \emt & 3 & 2 \\ 
   \bas{6} & 1 
}}
\put(265,60){\sqthree}
\put(265,42){\sqtwo}
\put(265,24){\sqone}
\put(265,6){\sqtwo}
\put(200,60){
$\longrightarrow$  \;
\tableau{ \bas{{10}} & {10} & {10} & {10} & 3 & 3 \\
\bas{{9}} & {9} & {9} & 1 \\
\bas{{8}} & {8} & 2 \\
\bas{{7}} & {7} & {7} & 3 & 2 \\
\bas{{6}} & 1
}}
\end{picture}
\]
\caption{Construction example for a pair $Y$, $\sigma$}
\label{fig-shape-alg-example}
\end{center}
\end{figure}

\begin{remark} \label{rem:LRS-algorithm}
The proof of Proposition \ref{!shape} provides us with an 
algorithm for constructing the desired SSK on a large basement
by successively filling the ``lowest'' row strip in the unfilled portion of the diagram for the set of columns containing the smallest-valued entries at each step,
as illustrated in Figure \ref{fig-shape-alg-example}.
An easy argument shows that starting with an SSK $Y$ as in the statement of the proposition,
if we have two compositions $\sigma$ and $\sigma'$, both permutations of $\delta$ with $\strongof{\sigma}=\strongof{\sigma'}$, 
then the respective constructed SSK $L$ and $L'$ will have respective shapes $\sigma/\tau$ and $\sigma'/\tau'$ with $\strongof{\tau}=\strongof{\tau'}$.
\end{remark}

%%%%%%%
\section{Littlewood-Richardson rule for Demazure atoms} \label{sec:LR-atoms}

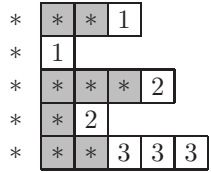
\begin{figure}[ht]
\begin{center}
\[
\begin{picture}(100,100)
\put(23,80){\sqtwo}
\put(23,44){\sqthree}
\put(23,26){\sqone}
\put(23,8){\sqtwo}
\put(0,80){\tableau{
\bas{{*}} & {*} & {*} & 1 \\
\bas{{*}} & 1 \\
\bas{{*}} & {*} & {*} & {*} & 2 \\
\bas{{*}} & {*} & 2 \\
\bas{{*}} & {*} & {*} & 3 & 3 & 3}}
\end{picture}
\]
\caption{An LRS with $n=5$ and column reading word $3231321$}
\label{fig-LRS}
\end{center}
\end{figure}

A \emph{Littlewood-Richardson skew skyline tableau} (LRS) of shape $\delta/\gamma$ is an 
SSK of shape $\delta/\gamma$ with large basement $b_i=2n-i+1$, where $n = \ell (\delta)=\ell (\gamma)$,
whose column reading word is a regular contre-lattice word.
Figure \ref{fig-LRS} shows an example of an LRS
with column reading word $3231321$, which is regular contre-lattice of content $(2,2,3)$.
We let $\text{LRS}(n)$ denote the set of LRS with entries in $[n]$.

We can now state our LR rule for the product of a Schur function and a Demazure atom.
\begin{theorem}  \label{thm:LR-Dzatoms}
In the expansion 
\begin{equation} \label{eq:LR-Dzatoms}
\dzatom_\gamma (X_n) \cdot s_\lambda (X_n)= 
\sum_\delta a^\delta_{\gamma\lambda} \dzatom_\delta (X_n),
\end{equation}
the coefficient $a^\delta_{\gamma\lambda}$ is the number of 
elements in $\text{LRS}(n)$ of shape $\delta/\gamma$ with content $\lambda ^{*}$.
\end{theorem}
%Before proving this theorem, 
%%% stuff moved!
%\bigskip
%%%
\begin{proof}%[Proof of Theorem \ref{thm:LR-Dzatoms}]
As with the proof of the classical LR rule for Schur functions~\cite{Fulton},
we recall the homomorphism $\psi : T \mapsto \xx^T$ from  the 
contretableau ring $R_n$, the graded algebra whose basis is $\text{SSYT}(n)$, onto the 
polynomial ring $\bZ[X] = \bZ[x_1,\ldots,x_n]$.
Under the bijection $\rho$ we may identify SSK with their corresponding CT.
The combinatorial formulas given in Equations  \eqref{eq:schurCT} and  \eqref{eq:atom} 
allow us to identify pre-images of Schur functions and Demazure atoms:
\begin{eqnarray}
S_\lambda &=& \sum_{V\in \text{CT}(n),\atop \text{shape}(V)=\lambda} V,  
\qquad\qquad \psi(S_\lambda)=s_\lambda(X_n) \\
A_\gamma &= &  \sum_{U\in \text{SSKI}(n),\atop \text{shape}(U)=\gamma} U, 
\qquad\qquad \psi(A_\gamma) =\dzatom_\gamma(X_n)
\end{eqnarray}
Under the homomorphism we then have $\psi(A_\gamma\cdot S_\lambda) = 
\dzatom_\gamma(X_n)\cdot s_\lambda(X_n)$.
The terms of $A_\gamma\cdot S_\lambda$ are the products of ordered pairs of 
CT $(U,V)$ where $\rho^{-1}(U)$ has shape $\gamma$.
The idea of the proof is then to exploit the the bijection $(U,V)\leftrightarrow(T,S)$ 
between ordered pairs $(U,V)$ of arbitrary CT and pairs $(T,S)$ of a CT $T = U\cdot V$ and a 
recording LR skew CT $S$, restricting to the case where the 
image $\rho^{-1}(U)$ has shape $\gamma$.

The bijection matches the CT $V$ (which here has shape $\lambda$) with a super CT of the same shape,  
and the pair is mapped to a biword $W$ using the RSK correspondence.  (See \cite[Chapter $7$]{ECII} for
a discussion of the RSK algorithm.
Note that for CT, the biword $W$ is in reverse lexicographic order.)
We then compute $(T,S) = U\gets W$.
That is, the lower row of the biword is \emph{inserted} into $U$ to obtain the pair $T = U\cdot V$ 
while the upper row of the biword is \emph{placed} into the corresponding skew Ferrers shape to 
obtain an LR skew CT $S$.
In the same way we can compute $(\rho^{-1}(T),L) = \rho^{-1}(U)\gets W$.
As we insert/place the biletters one-by-one, we can also track the images of the intermediate CT 
under the bijection $\rho^{-1}$, that is, inserting the bottom row of $W$ into the 
SSK $\rho^{-1}(U)$ using the insertion map described in~\cite{mason-1} and placing 
the upper row of $W$ in an SSK $L$, recording the location of the new cell.
Figure \ref{fig-LRatom-example} gives an example.
The resulting insertion SSK will of course be $\rho^{-1}(T)$, say of shape $\delta$.
It remains to show that (1) the resulting SSK $L$, 
when combined with the basement $b_i=2n-i+1$, is in fact an LRS, 
and (2) conversely, that any LRS $L$ of shape $\delta/\gamma$ and weight $\lambda ^{*}$ 
can be used to evacuate a biword $W$ from any SSK $\rho^{-1}(T)$ of shape $\delta$, 
leaving an SSK $\rho^{-1}(U)$ of shape $\gamma$, and such that the lower row of $W$ 
rectifies to a CT $V$ of shape $\lambda$ such that $T=U\cdot V$.

%\[ \begin{CD}
%   G  @>M>> H \\
%   @V\Delta_GVV   @VV\Delta V \\
%   G \otimes G  @>>M_G>  G
%   \end{CD} \]
\begin{figure}[ht]
\begin{center}
\[
\begin{array}{c}
\begin{pmatrix}
\;\tableau{  5 & 5 & 2 \\ 4 & 2 \\ 3 & 1 \\ 1 }  & \cdot &
\tableau{  5 & 3 & 2 \\ 4 & 2 \\ 3 & 1} \; \vspace{6pt}\\ U & & V
\end{pmatrix} 
\qquad = 
\vspace{12pt} \\
\begin{CD}
\begin{pmatrix} \;
\tableau{  5 & 5 & 2 \\ 4 & 2 \\ 3 & 1 \\ 1 } & \gets &
\begin{pmatrix}
3 & 3 & 3 & 2 & 2 & 1 & 1 \\
3 & 2 & 1 & 4 & 2 & 5 & 3
\end{pmatrix}  \vspace{6pt} \\  U & & W
\end{pmatrix} @>\rho^{-1}>>
\begin{pmatrix}
\tableau{
   \bas{1} & 1 & 1 \\
   \bas{2}  \\ 
   \bas{3} & 3 & 2 & 2 \\
   \bas{4} & 4 \\ 
   \bas{5} & 5 & 5  
} &\gets&
\begin{pmatrix}
3 & 3 & 3 & 2 & 2 & 1 & 1 \\
3 & 2 & 1 & 4 & 2 & 5 & 3
\end{pmatrix}   \vspace{6pt} \\  \rho^{-1}(U) & & W
\end{pmatrix} 
\\ @VVV   @VVV \\
\begin{pmatrix} \;
\tableau{5 & 5 & 5 & 3 & 2 \\ 4 & 4 & 2 & 2 \\ 3 & 3 & 1 \\ 2 & 1 \\ 1}
\;,  & 
\begin{picture}(100,0)
\put(0,0){\sqthree}
\put(0,-18){\sqtwo}
\put(0,-36){\sqtwo}
\put(0,-54){\sqone}
\put(0,0){\tableau{{*} & {*} & {*} & 3 & 3 \\ {*} & {*} & 3 & 2  \\ {*} & {*} & 1 \\ {*} & 2 \\ 1}}
\end{picture}
\; \vspace{6pt} \\ T & S
\end{pmatrix} @>\rho^{-1}>>
%\stackrel{\rho^{-1}}{\longleftrightarrow}
\begin{pmatrix}
\tableau{
\bas{1} & 1 & 1 & 1 \\
\bas{2} & 2 \\
\bas{3} & 3 & 3 & 2 & 2 & 2 \\
\bas{4} & 4 & 4 \\
\bas{5} & 5 & 5 & 5 & 3
}
\;, & 
\begin{picture}(120,20)
\put(23,0){\sqtwo}
\put(23,-36){\sqthree}
\put(23,-54){\sqone}
\put(23,-72){\sqtwo}
\put(0,0){\tableau{
\bas{{10}} & {10} & {10} & 1 \\
\bas{{9}} & 1 \\
\bas{{8}} & {8} & {8} & {8} & 3 & 3 \\
\bas{{7}} & {7} & 2 \\
\bas{{6}} & {6} & {6} 
}}
\end{picture}
 \vspace{6pt} \\ \rho^{-1}(T) & L
\end{pmatrix}  \\
\end{CD} 
\end{array}
\]
\caption{Correspondence example for a term of $\dzatom_{(2,0,3,1,2)}\cdot s_{(3,2,2)}$}
\label{fig-LRatom-example}
\end{center}
\end{figure}
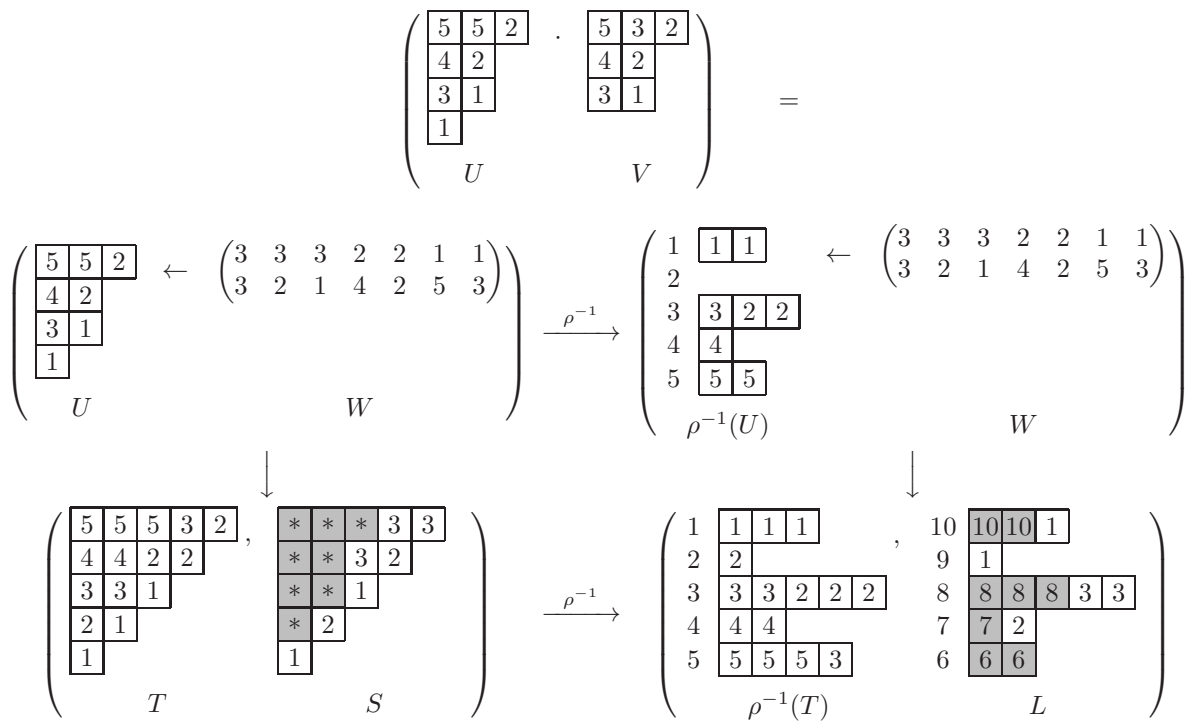

In the first direction, suppose we have constructed $L$ from $\rho^{-1}(U)\gets W$ as above.
We claim that $L$, including its basement, satisfies the conditions of an SSK.
Since we construct $\rho^{-1}(T)$ by adding successive row strips into $\rho^{-1}(U)$,
we likewise are constructing $L$ by adding successive row strips to the basement of shape $\gamma$.
This implies that the entries within each column of $L$ are distinct.
Since we place into $L$ the higher-numbered entries first, this forces the entries 
within each row of $L$ to be weakly decreasing left-to-right.

In the construction, suppose that after the addition of some particular cell,
row $j$ of the resulting SSK is strictly longer than row $i$ for some $i < j$.
As a consequence of the single box case of the Pieri rule,
it follows that at every following stage of the construction row $j$ must be strictly longer than row $i$,
and so $\delta_i < \delta_j$.
We consider triples in $L$.
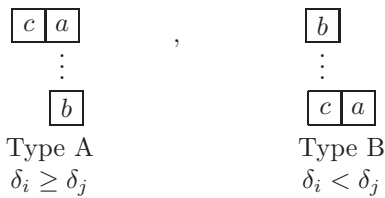
\[
\begin{array}{ccc} \vspace{6pt}
\tableau{   c & a   \\  & \bas{{\vdots}} } 
&  \qquad , \qquad  &
\tableau{   b   \\  \bas{{\vdots}}} 
 \\
 \tableau{ & b} & & \hspace*{.2in} \tableau{c & a} \\
\mbox{Type A} & & \hspace*{.2in} \mbox{Type B} \\
\delta_i \geq \delta_j & & \hspace*{.2in} \delta_i < \delta_j
\end{array}
\]
Suppose $L$ has a coinversion type A triple $\bigr((i,k),(j,k),(i,k-1)\bigr)$ with values $(a,b,c)$, as shown.
Since rows are weakly decreasing, this would imply $a<b \leq c$,
implying that the cell $(i,k)$ is not in the basement, and was filled after the cell $(j,k)$. 
But this would imply that just prior to adding cell $(i,k)$, row $j$ was longer than row $i$,
which in turn would imply that $\delta_i < \delta_j$,
contradicting that the three cells form a type A triple.  
Thus $L$ can have no type A coinversion triples.

Suppose $L$ has a coinversion type B triple $\bigr((j,k+1),(i,k),(j,k)\bigr)$ with values $(a,b,c)$, as shown.
Since rows are weakly decreasing, this would imply $a \leq b<c$.
But that would then imply that cell $(i,k)$ is not in the basement, and was added 
before cell $(j,k+1)$, a violation of the Pieri rule.
Thus $L$ can have no type B coinversion triples,
and so $L$ is a valid SSK.

To see that  $\text{colword}(L)$ is a regular contre-lattice word,  note that 
within each column of $L$, the set of values (excluding the basement) in that column is the 
same as the set of values in the corresponding column of the LR skew CT $S$.
We may consider $S$ to be an SSK with basement $b_i=2n-i+1$ of shape $\partitionof{\gamma}$.
Since   $\text{colword}(S)$ is a regular contre-lattice word,
by Proposition \ref{prop:r-contre}, $S$ is loosely contre-lattice, where the maximum entry in
$\text{colword}(S)$ is $\ell(\lambda)$.
Since $L$ has the same column sets as $S$, $L$ is also loosely contre-lattice,
and by Proposition \ref{prop:r-contre} $L$ is contre-lattice, i.e.
$\text{colword}(L)$ is a regular contre-lattice word, and so $L$ is an LRS as claimed.

\medskip
For the converse direction, assume $L$ is an LRS of shape  $\delta/\gamma$ and 
weight $\lambda ^{*}$, and that $\rho^{-1}(T)$ is any SSK of shape $\delta$.
We show that we can use $L$ to evacuate a CT from $\rho^{-1}(T)$ as desired.
In the process we construct a biword $W$.
We know from  Proposition \ref{prop:r-LRS-reduction} that the rightmost least entry 
in $L$ having entry value 1, call it  $x_1$, appears at the end of the last row in $L$ of 
some particular length, which by the Pieri rule implies that we can evacuate the corresponding 
cell of $\rho^{-1}(T)$, obtaining a value $v$ and leaving a new SSK $\rho^{-1}(T')$ of the 
same shape as $L' = L-x_1$.  
We record  ${1\choose v}$ as the last biletter of $W$.
Since by Proposition \ref{prop:r-any-reduction} the remaining SSK $L'$ is also contre-lattice, we can repeat the process, 
constructing the biletters of $W$ in reverse order, until all cells of $L$ and their 
corresponding cells in $\rho^{-1}(T)$ have been processed, leaving us with a remaining 
SSK $\rho^{-1}(U)$ of shape $\gamma$ and a biword $W$ whose upper row has weight $\lambda ^{*}$.
$W$ and $\rho^{-1}(U)$  in turn correspond to a pair of CT $(V,H)$ of the same shape, 
where clearly $T = U\cdot V$.
To see that in fact $H$ is the super CT of shape $\lambda$,
consider the parallel step-by-step construction using $S$ to evacuate $T$, 
where $S$ is the LR skew CT of shape $\partitionof{\delta}/\partitionof{\gamma}$ obtained by 
sorting the columns of $L$ (including the basement), as illustrated in Figure \ref{fig-LRatom-example}.
Evacuating a cell of $\rho^{-1}(T)$ corresponding to a cell $x_k$ in $L$ (containing the $k$th entry of 
value $x$ in column reading order)
corresponds under the bijection $\rho$ to evacuating a cell of $T$ corresponding to the 
cell $x_k$ in $S$, producing the same biword $W$.
Since $S$ rectifies to the super CT of shape $\lambda$, $V$ also has shape $\lambda$.
\end{proof}

%%%%%%%
\section{Littlewood-Richardson rule for quasisymmetric Schur functions}
\label{LRQS}

Consider an $\text{SSK}$ with basement $b_i=i$.  
It is easy to see that if $\gamma _k>0$, then in any such SSK $T$ of shape $\gamma$, the cell of $\gamma$ in column one and row $k$ must contain the number $k$.  
If we consider SSK with an arbitrary increasing basement \[ 1\leq b_1 < b_2 < \cdots < b_n, \]
where the $b_i$ are not necessarily consecutive, 
then we can identify $T$ with a unique SSK $\hat{T}$ of shape $\strongof{\gamma}$ obtained by removing the rows of zero length.
For such SSK of composition shape, the basement becomes superfluous.
This motivates the following definition.

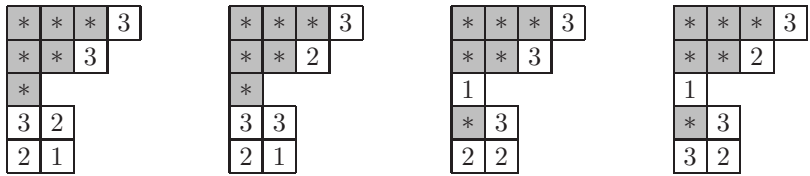
\begin{figure}[htbp]
\begin{center}
\[
\begin{picture}(440,100)
\put(0,80){\sqthree}
\put(0,62){\sqtwo}
\put(0,44){\sqone}
\put(0,80){\tableau{{*} & {*} & {*} & 3 \\
{*} & {*} & 3 \\
{*} \\
3 & 2 \\
2 & 1 }}
\put(120,80){\sqthree}
\put(120,62){\sqother}
\put(120,44){\sqone}
\put(120,80){\tableau{{*} & {*} & {*} & 3 \\
{*} & {*} & 2 \\
{*} \\
3 & 3 \\
2 & 1 }}
\put(240,80){\sqthree}
\put(240,62){\sqother}
\put(240,26){\sqone}
\put(240,80){\tableau{{*} & {*} & {*} & 3 \\
{*} & {*} & 3 \\
1 \\
{*} & 3 \\
2 & 2 }}
\put(360,80){\sqthree}
\put(360,62){\sqtwo}
\put(360,26){\sqone}
\put(360,80){\tableau{{*} & {*} & {*} & 3 \\
{*} & {*} & 2 \\
1 \\
{*} & 3 \\
3 & 2 }}
\end{picture}
\]
\caption{The four distinct LRC of shape $(4,3,1,2,2)/(3,2,1)$ and content $(1,2,3)$}
\label{fig-LRC}
\end{center}
\end{figure}

%Let $\beta$ be a composition and $\tau$ be a weak composition with $\ell(\beta) = \ell (\tau) =m<n$, and assume $\beta _i > \tau_i$ for $1\le i \le m$.
We define a \emph{semistandard composition tableau} ($\text{SSC}$) of shape
$\beta$ (a composition) to be a filling of the diagram $\beta$ which is 
strictly increasing down the first column, weakly decreasing rightward along
each row, and where every triple is an inversion triple.
Since by definition the \emph {QS function} $\qschur_\alpha$ is the sum of Demazure atoms $\dzatom _{\gamma}$, 
over all $\gamma$ with $\strongof {\gamma}=\alpha$, (\ref{eq:atom}) implies that \cite{HLMvW}
\begin{equation} 
\label{eq:qschur-combo}
\qschur_\alpha(X_n) = \sum_{T\in SSC(n),\atop shape(T)=\alpha} \xx^T,
\end{equation}
where $\text{SSC}(n)$ is the set of all $\text{SSC}$ with entries in $[n]$.  
These QS functions also satisfy an LR rule.  To state it we need to define
an analog of LRS.
Let $L_1$ and $L_2$ be elements of LRS($n$), where $L_1$ has shape $\delta/\sigma$ and $L_2$ has shape $\gamma/\tau$.
We declare $L_1$ and $L_2$ to be equivalent if
\begin{enumerate}
\item  $L_1$ and $L_2$ have the same set of non-basement entries in each column.
\item  $\strongof{\delta}=\strongof{\gamma}$, say $\strongof{\delta}=\strongof{\gamma}=\beta$.
%\item  $\strongof{\sigma}=\strongof{\tau}$, say $\strongof{\sigma}=\strongof{\tau}=\alpha$.
\end{enumerate}
We define a \emph{Littlewood-Richardson skew composition tableau} (LRC) to be an equivalence class of $\text{LRS}(n)$, and the collection of such equivalence classes we denote $\text{LRC}(n)$.
Each equivalence class determines a sequence of column sets and a pair of compositions $\beta$ and $\alpha$ which are the underlying compositions of the overall shape and basement shape respectively of the elements of the equivalence class.
We shall define the \emph{shape} of the LRC to be this pair of compositions
and by abuse of notation we shall denote the shape by $\beta/\alpha$.
(In view of Remark \ref{rem:LRS-algorithm}, the shapes of the respective basements of the elements of a given $\text{LRC}$
equivalence class all have the same underlying partition $\alpha$, and hence the shape $\beta /\alpha$ is
well-defined.)
We can represent an LRC diagrammatically.
In Figure \ref{fig-LRC} we exhibit the four LRC of shape $(4,3,1,2,2)/(3,2,1)$ and content $(1,2,3)$.

We can now state the LR rule for the product of a QS function and a Schur function.

\begin{theorem}  \label{thm:LR-qschur}
In the expansion 
\begin{equation} \label{eq:LR-qschur}
\qschur_\alpha(X_n) \cdot s_\lambda(X_n) = \sum_\beta C^ \beta_{\alpha\lambda} \qschur_\beta(X_n),
\end{equation}
the coefficient $C^ \beta_{\alpha\lambda}$ is the number of elements in $\text{LRC}(n)$ of 
shape $\beta/\alpha$ with content 
$\lambda ^{*}$.
\end{theorem}

\begin{proof}
We make use of \eqref{eq:qschur-combo} 
and \eqref{eq:schurCT}.
The $\text{SSC}$ are trivially in bijection with the SSKI, that is, the SSK with $b_i=i$, 
hence the mapping $\rho:\text{SSKI}(n)\to \text{CT}(n)$ can be viewed as a
bijection between $\text{SSC}$ $D$ and CT $\rho(D)$ whose 
columns are just the respective sorted column sets of $D$.
In view of the proof of Theorem \ref{thm:LR-Dzatoms},
it suffices to provide a 
bijection $(U,V)\leftrightarrow(T,S)$ between pairs $(U,V)$ of 
$\text{CT}$, where $\rho ^{-1}(U)$ is an $\text{SSC}$ of shape $\alpha$ and $V$ has 
shape $\lambda$, and 
pairs $(T,S)$, where $T$ is the CT $T=U\cdot V$, with $\rho ^{-1}(T)$ an 
$\text{SSC}$ of shape $\beta$, and $S$ is an LRC of weight $\lambda ^{*}$ and 
shape $\beta/\alpha$.

We make use of the bijection $(U,V)\leftrightarrow(T,L)$ as constructed in the proof of  Theorem  \ref{thm:LR-Dzatoms}.
Suppose here, as in the proof, that $(U,V)$ is a pair of CT with entries in
$[n]$.
Under the bijection $\rho^{-1}:\text{CT}(n)\to \text{SSKI}(n)$ 
we can map $U$ and $T=U\cdot V$ respectively to SSKI with $n$ rows, say $\rho^{-1}(U)$ of shape $\gamma$ and $\rho^{-1}(T)$ of shape $\delta$.
Under the bijection from the proof of Theorem \ref{thm:LR-Dzatoms}, $\rho^{-1}(T)$ is 
paired with an LRS $L$ of weight $\lambda ^{*}$ and shape $\delta/\gamma$.
Thus the pair $(U,V)$ determines a unique pair $(T,S)$ where $S$ is the LRC of shape $\strongof{\delta}/\strongof{\gamma}$ that is the equivalence class of $L$.

Conversely, suppose we have a pair $(T,S)$ where $\rho^{-1}(T)$, viewed as an SSC,  has shape $\beta$ and $S$ is an element of $\text{LRC}(n)$ of weight 
$\reverse{\lambda}$ and of shape $\beta/\alpha$.
Then under the bijection $\rho^{-1}:\text{CT}(n)\to \text{SSKI}(n)$,  
$\rho^{-1}(T)$ is an SSKI of shape $\delta$, where $\ell(\delta)=n$ and $\strongof{\delta}=\beta$. 
%We claim that there is a unique LRS $L$ of the LRC (equivalence class) $S$ of overall shape $\delta$.
Now by Proposition \ref{!shape} the existence of $S$ implies that there is a unique LRS $L$ of shape $\delta/\gamma$ for some $\gamma $ with $\partitionof{\gamma}= \partitionof{\alpha}$ and having the same column sets of entries as the elements of $S$.
As mentioned in Remark \ref{rem:LRS-algorithm}, 
the construction in the proof of Proposition \ref{!shape} implies that $\strongof{\gamma}=\alpha$, that is, $L$ is in fact an element of $S$.
Thus the pair $(T,S)$ determines a unique pair $(T,L)$.
Furthermore, $L$ has weight $\reverse{\lambda}$.
Under the bijection from the proof of Theorem  \ref{thm:LR-Dzatoms}, $(T,L)$ is paired with a 
pair of CT $(U,V)$ where $T=U\cdot V$, $V$ has shape $\lambda$, and $\rho^{-1}(U)$ has 
shape $\gamma$, which implies that $\rho^{-1}(U)$, when viewed as an SSC, has shape $\strongof{\gamma}=\alpha$, as desired.
\end{proof}

%%%%%%%
\section{Littlewood-Richardson rule for Demazure characters} \label{sec:LR-chars}

A \emph{Littlewood-Richardson skew key} (LRK) of shape $\delta/\gamma$ is an 
SSK of shape $\delta/\gamma$ with basement $b_i=n+i$, where $n = \ell (\delta)=\ell (\gamma)$ and
whose column reading word is a regular contre-lattice word.
We let $\text{LRK}(n)$ denote the set of $\text{LRK}$ with entries in $[n]$.
Figure \ref{fig-LRK} provides an example of an LRK of shape $(5,1,3,2,4)/(2,0,1,2,3)$ and
$\text{colword}=3323121$.  

We can now state our LR rule for the product of a Schur function 
and a Demazure character.

\begin{figure}[ht]
\begin{center}
\[
\begin{picture}(100,100)
\put(23,80){\sqtwo}
\put(23,44){\sqone}
\put(23,26){\sqtwo}
\put(23,8){\sqthree}
\put(0,80){\tableau{
\bas{{6}} & {6} &{6} & 3 & 3 & 3 \\
\bas{{7}} & 1 \\
\bas{{8}} & {8} & 2 & 1 \\
\bas{{9}} &{9} & {9} \\
\bas{{10}} & {10} & {10} & {10} & 2}}
\end{picture}
\]
\caption{An LRK with $n=5$ and column reading word $3323121$}
\label{fig-LRK}
\end{center}
\end{figure}

\begin{theorem}  \label{thm:LR-Dzchars}
In the expansion 
\begin{equation} \label{eq:LR-Dzchars}
\dzchar_\gamma(X_n) \cdot s_\lambda(X_n) = \sum_\delta b^\delta_{\gamma\lambda} \dzchar_\delta(X_n),
\end{equation}
the coefficient $b^\delta_{\gamma\lambda}$ is the number of elements in 
$\text{LRK}(n)$ of shape $\delta^*/\gamma^*$ with content $\reverse{\lambda}$.
\end{theorem}

\begin{proof}[Proof of Theorem \ref{thm:LR-Dzchars}]

Recall \cite{LaSc}, \cite{mason-2} that the Demazure characters can be obtained from the Demazure atoms:\begin{equation}
\dzchar_\gamma = \sum_{\beta \geq \gamma^*} \dzatom_\beta,
\end{equation}
where the sum is over all compositions $\beta$ which are weakly above $\gamma ^*$ in the Bruhat order.
(Given a weak composition $\gamma$, let $\pi(\gamma)$ be the permutation of minimal length which 
arranges the parts of $\gamma$ into nonincreasing order.
Then we define $\beta \ge \alpha$ if and only if $\pi (\beta)\le \pi(\alpha)$ in the usual (strong) Bruhat order on permutations.)
We substitute the formula \eqref{eq:LR-Dzatoms} for the multiplication of a Demazure atom 
and a Schur function to obtain the following formula for the left hand side of \eqref{eq:LR-Dzchars}: 
$$
\dzchar_{\gamma} \cdot s_\lambda %= \left( \sum_{\beta \geq \reverse{\gamma}} \dzatom_\beta \right) \cdot s_\lambda 
= \sum_{\beta \geq \reverse{\gamma}} \dzatom_\beta \cdot s_\lambda = \sum_{\beta \geq \reverse{\gamma}} 
\sum_\delta a^\delta_{\beta\lambda} \dzatom_\delta,
$$
where $a^\delta_{\beta\lambda}$ is the number of elements of $\text{LRS}(n)$ of shape $\delta / \beta$ 
with content $\reverse{\lambda}$.   To prove that
\begin{equation} \label{eq:dzchar-2}
\sum_\alpha b^\alpha_{\gamma\lambda} \dzchar_\alpha = \sum_{\beta \geq \reverse{\gamma}} 
\sum_\delta a^\delta_{\beta\lambda} \dzatom_\delta, 
\end{equation}
we further expand the left hand side of \eqref{eq:dzchar-2} to see that our theorem is equivalent to the identity
\begin{equation}{\label{dem-atom-equiv}}
\sum_\alpha b^\alpha_{\gamma\lambda} \sum_{\delta \geq \reverse{\alpha}} \dzatom_\delta  
= \sum_{\delta} \sum_{\beta \geq \reverse{\gamma}} a^\delta_{\beta\lambda} \dzatom_\delta.
\end{equation}
Each coefficient $b^\alpha_{\gamma\lambda}$ appearing on the left hand side of (\ref{dem-atom-equiv}) is the coefficient of every Demazure atom $\dzatom_\delta$ such that $\delta \geq \reverse{\alpha}$.  
Since the Demazure atoms are linearly independent, comparing the coefficients of
$\mathcal{A}_{\delta}$ on both sides of (\ref{dem-atom-equiv}) reduces our identity to
\begin{equation}{\label{dmcprod}}
 \sum_{\delta\geq\reverse{\alpha}\supseteq\reverse{\gamma}} b^\alpha_{\gamma\lambda} =
 \sum_{\delta\supseteq\beta\geq \reverse{\gamma}} a^\delta_{\beta\lambda}
\end{equation}
for fixed $\delta$ and $\gamma$.
It therefore suffices to fix $\delta$ and $\gamma$ and find a bijection between the set $K$
of all LRK of shape $\reverse{\alpha}/\reverse{\gamma}$ with content $\reverse{\lambda}$ where $\reverse{\alpha} \le \delta$ in Bruhat order 
and the set $L$ of all LRS of shape $\delta/\beta$ with content $\reverse{\lambda}$ where $\beta\geq\reverse{\gamma}$ in Bruhat order.

We begin with the forward direction of the map $\phi: K \mapsto L$.  
Let $K$ be an LRK with content  $\reverse{\lambda}$ and shape $\reverse{\alpha}/\reverse{\gamma}$, and assume $\reverse{\alpha} \le \delta$ in Bruhat order.  
By Proposition \ref{!shape} there exists a unique LRS $L$ of shape $\delta/\beta$ for some $\beta$ a permutation of $\gamma$ and having the same column sets as $K$.  Map the LRK $K$ to this LRS $L$.
To show that the map takes $K$ into the appropriate set, 
we must prove that $\beta\ge\reverse{\gamma}$  in Bruhat order.

To see this, let $\gamma _K$ be the overall shape of $K$ and $\gamma_L$ be the overall
shape of $L$ and apply the following 
iterative argument.  The overall shape of $L$ is weakly higher than the reverse of the shape of $K$ by construction, 
so $\reverse{\gamma_K} \le \gamma_L$. 
In the construction given by Proposition \ref{!shape}, 
consider the first entry in $K$ that is mapped to $L$.  
This entry is mapped to a row $(r_1)_L$ of $L$ weakly higher than the reverse of the row $(r_1)_K$ of $K$ from which it is removed since  the largest part of $\gamma_L$ appears 
before the largest part of $\reverse{\gamma_K}$.  
Subtract one from the $(r_1)_L$ part of $\gamma_L$ and the $(r_1)_K$ part of $\reverse{\gamma_K}$ to obtain new compositions $\reverse{\gamma_K} \le \gamma_L$. 
Repeat this procedure until there are no remaining non-basement entries in $K$.  
The resulting compositions $\reverse{\gamma_K}$ and $\gamma_L$ are the shapes of the respective basements and satisfy 
$\reverse{\gamma_K} \le \gamma_L$.  Therefore the basement of $L$ is indeed higher in Bruhat order than the 
reverse basement of $K$.  
(See Figure \ref{fig-bruhat} for an example.)

\begin{figure}[htbp]
\begin{center}
\[
\begin{array}{ccccc}
K & \gamma_K^* & & L & \gamma_L  
\vspace{6pt}\\ \hline
\\
% Row 1 of array
\begin{picture}(100,0)
\put(23,0){\sqone}
\put(23,-36){\sqtwo}
\put(23,-54){\sqone}
\put(23,-72){\sqthree}
\put(0,0){\tableau{
\bas{{6}} & {6} & 2 \\
\bas{{7}} \\
\bas{{8}} & {8} & {8}
& 3 & 3 \\
\bas{{9}} & {9} & 1 \\
\bas{{10}} & {10} & {10} & {10}}}
\end{picture}
& (3,2,4,0,2)   
& \quad \text{ } \quad &
\tableau{
\bas{10} & {} & {} & {} \\
\bas{9} & {} & {} & {} & {} \\
\bas{8} & {} & {} \\
\bas{7} \\
\bas{6} & {} & {}  }
 & (3,4,2,0,2) 
\vspace{20pt} \\
% Row 2 of array
\begin{picture}(100,0)
\put(23,0){\sqone}
\put(23,-36){\sqtwo}
\put(23,-54){\sqone}
\put(23,-72){\sqthree}
\put(0,0){\tableau{
\bas{{6}} & {6} & 2 \\
\bas{{7}} \\
\bas{{8}} & {8} & {8}
& 3 & 3 \\
\bas{{9}} & {9} \\
\bas{{10}} & {10} & {10} & {10}}}
\end{picture}
 & (3,1,4,0,2)  & &
\tableau{
\bas{10} & {} & {} & {} \\
\bas{9} & {} & {} & {} & {} \\
\bas{8} & {} & {} \\
\bas{7} \\
\bas{6} & {} & {1}  }
 & (3,4,2,0,1)  
\vspace{20pt} \\
% Row 3 of array
\begin{picture}(100,0)
\put(23,0){\sqone}
\put(23,-36){\sqtwo}
\put(23,-54){\sqone}
\put(23,-72){\sqthree}
\put(0,0){\tableau{
\bas{{6}} & {6}  \\
\bas{{7}} \\
\bas{{8}} & {8} & {8}
& 3 & 3 \\
\bas{{9}} & {9}  \\
\bas{{10}} & {10} & {10} & {10}}}
\end{picture}
 & (3,1,4,0,1)  & &
\tableau{
\bas{10} & {} & {} & {} \\
\bas{9} & {} & {} & {} & {} \\
\bas{8} & {} & {2} \\
\bas{7} \\
\bas{6} & {} & {1}  }
& (3,4,1,0,1)  
\vspace{20pt} \\
% Row 4 of array
\begin{picture}(100,0)
\put(23,0){\sqone}
\put(23,-36){\sqtwo}
\put(23,-54){\sqone}
\put(23,-72){\sqthree}
\put(0,0){\tableau{
\bas{{6}} & {6}  \\
\bas{{7}} \\
\bas{{8}} & {8} & {8}
& 3  \\
\bas{{9}} & {9}  \\
\bas{{10}} & {10} & {10} & {10}}}
\end{picture}
 & (3,1,3,0,1)  & &
\tableau{
\bas{10} & {} & {} & {} \\
\bas{9} & {} & {} & {} & {3} \\
\bas{8} & {} & {2} \\
\bas{7} \\
\bas{6} & {} & {1}  }
 & (3,3,1,0,1)  
\vspace{20pt} \\
% Row 6 of array
\begin{picture}(100,0)
\put(23,0){\sqone}
\put(23,-36){\sqtwo}
\put(23,-54){\sqone}
\put(23,-72){\sqthree}
\put(0,0){\tableau{
\bas{{6}} & {6}  \\
\bas{{7}} \\
\bas{{8}} & {8} & {8} \\
\bas{{9}} & {9}  \\
\bas{{10}} & {10} & {10} & {10}}}
\end{picture}
 & (3,1,2,0,1) & &
\tableau{
\bas{10} & {} & {} & {} \\
\bas{9} & {} & {} & {3} & {3} \\
\bas{8} & {} & {2} \\
\bas{7} \\
\bas{6} & {} & {1}  }
 & (3,2,1,0,1) 
 \vspace{12pt} \\
 & & &
\begin{picture}(100,100)
\put(23,80){\sqthree}
\put(23,62){\sqtwo}
\put(23,44){\sqone}
\put(23,8){\sqone}
\put(0,80){\tableau{
\bas{{10}} & {{10}} & {{10}} & {{10}} \\
\bas{{9}} & {{9}} & {{9}} & {3} & {3} \\
\bas{{8}} & {{8}} & {2} \\
\bas{{7}} \\
\bas{{6}} & {{6}} & {1}  }}
\end{picture}
\end{array}
\]
\caption{Constructive comparison of the basements of $K$ and $L$}
\label{fig-bruhat}
\end{center}
\end{figure}

We have now shown that $L$ is an $\text{LRS}$ in the desired set.
%whose basement is weakly higher in Bruhat order than the reverse of the basement of $K$. 
Proposition \ref{!shape} shows that $L$ is unique, therefore the map $\phi$ is injective.  
We describe the inverse of the map $\phi$ to prove that the map is surjective.  
Consider an arbitrary $\text{LRS}$ $L$ of shape $\delta/\beta$ and 
content $\reverse{\lambda}$,  % Let $\delta$ denote the shape of $L$.  
and let $\gamma$ be a rearrangement of $\beta$ such that $\reverse{\gamma}\leq\beta$ in Bruhat order.
We need to map $L$ back to an $\text{LRK}$ of shape $\reverse{\alpha}/\reverse{\gamma}$ for some $\reverse{\alpha}\leq\delta$ in Bruhat order and having the same column sets as $L$.
%a given basement whose shape $\reverse{\alpha}$ is a rearrangement of $\beta$ and whose reverse is lower in Bruhat order than $\beta$.  
(Note that all $\text{LRK}$ in the pre-image have the same fixed basement of shape $\reverse{\gamma}$.)  
Let $K^0$ be the basement diagram of type $b_i=n+i$ and of shape $\reverse{\gamma}$.
This is the basement on which the $\text{LRK}$ will be built. 
Begin with the leftmost column of the $\text{LRS}$ $L$ 
and the largest non-basement entry in this column.  
Place this entry in the highest available row of this column in $K^0$, i.e. in an empty cell not part of the basement such that the entry to its left is non-empty and greater than or equal to our insertion entry.
Call the resulting SSK $K^1$.  
Repeat with the second largest entry in the leftmost column of $L$ to create $K^2$.  
Continue this procedure until all of the non-basement entries in the leftmost column of $L$ have been placed 
into the skyline diagram.  
Repeat with each column of $L$ from left to right until all of the non-basement entries of $L$ have 
been inserted into the SSK $K$.  
We must prove that $K$ is indeed an LRK, say of shape $\reverse{\alpha}/\reverse{\gamma}$,
 and that $\reverse{\alpha}\leq\delta$ in Bruhat order.
%such that the reverse of the shape of $K$ is below the shape $\delta$ in Bruhat order.  

%To see that the column reading word of $K$ is regular contre-lattice, consider the appearance of the non-basement entries in a given column of $K$.  Note that the entries must appear in the same column as in $L$, so if the column reading word of $K$ is not regular contre-lattice, then there exist entries $i$ and $i+1$ in the same column such that $i$ appears in a higher row of $K$ than $i+1$ and $i$ and $i+1$ appear an equal number of times to the right of this column.   But $i+1$ is inserted first, into the highest possible column, so this cannot happen.  Therefore the column reading word of $K$ must be regular contre-lattice.

The rows of $K$ are weakly decreasing by construction, so we must check that the triple conditions are satisfied.  We consider triples in $K$.
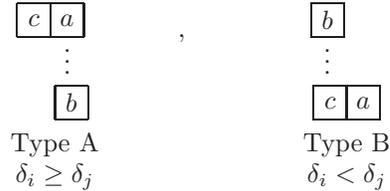
\[
\begin{array}{ccc} \vspace{6pt}
\tableau{   c & a   \\  & \bas{{\vdots}} } 
&  \qquad , \qquad  &
\tableau{   b   \\  \bas{{\vdots}}} 
 \\
 \tableau{ & b} & & \hspace*{.2in} \tableau{c & a} \\
\mbox{Type A} & & \hspace*{.2in} \mbox{Type B} \\
\delta_i \geq \delta_j & & \hspace*{.2in} \delta_i < \delta_j
\end{array}
\]

Suppose $K$ has a coinversion type A triple $((i,k),(j,k),(i,k-1))$ with values $(a,b,c)$ as shown.  Since the rows of $K$ are weakly decreasing, this would imply that $a < b \le c$ and therefore that the cell $(j,k)$ is not in the basement and was filled before the cell $(i,k)$.  But since the cell $(i,k-1)$ was filled before $(i,k)$ was filled, the entry $b$ would have been inserted into the cell $(i,k)$, a contradiction.  Thus $K$ can have no type A coinversion triples.

Suppose $K$ has a coinversion type $B$ triple $((j,k+1),(i,k),(j,k))$ with values $(a,b,c)$, as shown.  Since the rows of $K$ are weakly decreasing, this would imply $a \le b < c$.  That would then imply that the cell $(j,k+1)$ is not in the basement, and was added after the cell $(i,k+1)$, for otherwise the entry $a$ would be inserted into the cell $(i,k+1)$.  Therefore the entry in cell $(i,k+1)$ is greater than the entry in $(j,k+1)$ and will be filled first.  Continuing in this manner implies that $\delta_i \ge \delta_j$, contradicting the assumption that the three cells form a type B triple.  Thus $K$ can have no type B coinversion triples.

We invoke Proposition \ref{prop:r-contre} to see that the the diagram $K$ is contre-lattice, 
since $L$ is contre-lattice and the map from $L$ to $K$ preserves the column sets of the diagrams.  
We claim furthermore that $\reverse{\alpha}\leq\delta$ in Bruhat order.
To begin with, we have by assumption $\reverse{\gamma}\leq\beta$, where $\reverse{\gamma}$ is the shape of $K^0$, the basement of $K$, and $\beta$ is the shape of the basement of $L$.
As the first non-basement entry of $L$ is mapped to produce $K^1$, it will appear in a weakly higher row of $K^1$ than its appearance in $L$ by construction.  
The resulting shape of $K^1$ will therefore remain weakly lower in Bruhat order than the 
union of the basement of $L$ and this entry.  
Iterating this argument implies that the overall shape $\reverse{\alpha}$ of $K$ is 
weakly lower in Bruhat order than the overall shape $\delta$ of $L$. 
\end{proof}

\subsection{Recovering the classical Littlewood-Richardson rule}

Every Schur function is a Demazure character; in particular
$s_{\mu}(X_n)=\kappa_{\mu^*}(X_n)$.    Theorem \ref{thm:LR-Dzchars} is
therefore a generalization of the classical Littlewood-Richardson rule.
Consider the product
\begin{eqnarray*}
s_{\mu}(X_n) \cdot s_{\lambda}(X_n) & = & \kappa_{\mu^*}(X_n) \cdot
s_{\lambda}(X_n) \\
 & = & \sum_{\delta} b_{\mu^* \lambda}^{\delta} \kappa_{\delta}(X_n).
\end{eqnarray*}
We claim that $\sum_{\delta} b_{\mu^* \lambda}^{\delta}
\kappa_{\delta}(X_n) = \sum_{\nu} c_{\mu \lambda}^{\nu} s_{\nu}(X_n)$,
where $c_{\mu \lambda}^{\nu}$ is the number of LR skew CT with shape
$\nu/\mu$ and content $\lambda ^*$.  To see this, let $L$ be an arbitrary
element in $LRK(n)$ of shape $\delta^* / \mu$ and content $\lambda^*$.  The
basement of $L$ is the partition $\mu$.  If the shape of $\delta^*$ is not
a partition, then consider two rows $i$ and $j$ of $\delta^*$ such that
$i<j$ but row $j$ is strictly longer than row $i$.  Let $C$ be the column
containing the rightmost entry of the  basement in row $j$.  This entry,
together with the entry immediately to its right and the entry in row $i$
of column $C$ form a type $B$ coinversion triple.  Therefore the shape
$\delta^*$ must be a partition, and so $\kappa_{\delta}(X_n)$ is the Schur
function $s_{\delta^*}(X_n)$.

We already know that the row entries of $L$ weakly decrease left-to-right,
and, since $\delta ^*$ and $\mu$ are partitions, all inversion triples must
be of type $A$. Consequently our non-basement column entries decrease
top-to-bottom. Thus, $L$ is a skew CT. Lastly, note that $\text{colword}(L)$ is
regular contre-lattice if and only if $L$ is furthermore a LR skew CT by Proposition \ref{prop:LRcontretabs}. Therefore Theorem
\ref{thm:LR-Dzchars} reduces to the classical Littlewood-Richardson rule
whenever $\kappa_{\gamma}(X_n)$ is a Schur function.

\section{Acknowledgments}
The first author was supported in part by  NSF grant DMS-0553619 and DMS-0901467. The third author was supported in part by NSF postdoctoral research fellowship DMS-0603351. The second and fourth authors were supported in part by the National Sciences and Engineering Research Council of Canada. The authors would like to thank the Banff International Research Station and the Centre de Recherches Math\'{e}matiques, where some of the research took place.  The authors would also like to thank the referee for insightful comments and suggestions.

\end{document}